\documentclass[twosided]{amsart}
\usepackage{amsfonts,amssymb}
\newtheorem{thm}{Theorem}[section]
\newtheorem{cor}[thm]{Corollary}
\newtheorem{lem}[thm]{Lemma}

\newtheorem{defn}[thm]{Definition}

\begin{document}

\title{ $H^*$-algebras and Quantization of para-Hermitian spaces}%
\author{Gerrit van DIJK, Michael PEVZNER}%

\address{G.v.D.: Mathematish Instituut, Universiteit Leiden, PO Box 9512, NL-2300 RA
Leiden, Nederland, M.P.: Laboratoire de Math\'ematiques, UMR CNRS
6056,
 Universit\'e de Reims, Campus Moulin de la Housse BP 1039,
  F-51687, Reims, France.}%
\email{dijk{@}math.leidenuniv.nl, pevzner@univ-reims.fr}%


\begin{abstract}
In the present note we describe a family of $H^*$-algebra
structures on the set $ L^2(X)$ of square integrable functions on a
rank-one para-Hermitian symmetric space $X$.
\end{abstract}
\maketitle
\section*{Introduction}
Let $X$ be a para-Hermitian symmetric space of rank one. It is well-known that  $X$ is isomorphic (up to a covering) to the quotient space $SL(n,\mathbb
R)/GL(n-1,\mathbb R)$, see \cite{[KK]} for more details. We shall
thus assume throughout this note that $X=G/H$, where $G=SL(n,\mathbb
R)$ and $H =GL(n-1,\mathbb R)$.

The space $X$ allows the definition of a covariant symbolic calculus
that generalizes the so-called convolution-first calculus on
$\mathbb R^2$, see (\cite{[DM],[PU],[UU]}) for instance. Such a
calculus, or quantization map $\mathrm{Op}_\sigma$, from the set of
functions on $X$, called symbols, onto the set of linear operators
acting on the representation space of the maximal degenerate series
$\pi_{-\frac n2+i\sigma}$ of the group $G$, induces a
non-commutative algebra structure on the set of symbols, that we
suppose to be square integrable. On the other hand, the taking of
the adjoint of an operator in such a calculus defines an involution
on symbols. It turns out that these two data give rise to a
$H^*$-algebra structure on $L^2(X)$.

According to the general theory, (\cite{[A],[L],[N]}), every
$H^*$-algebra is the direct orthogonal sum of its closed minimal
two-sided ideals which are simple $H^*$-algebras. The main result of
this note is the explicit description of such a decomposition for
the Hilbert algebra $L^2(X)$ and its commutative subalgebra of
$SO(n,\mathbb R)$-invariants.
\section{Definitions and basic facts}
\subsection{$H^*$-algebras}
\begin{defn}
A set $R$ is called a $H^*$-algebra (or Hilbert algebra) if
\begin{enumerate}
  \item $R$ is a Banach algebra with involution;
  \item $R$ is a Hilbert space;
  \item the norm on the algebra $R$ coincides with the norm on the
  Hilbert space $R$;
  \item For all $x,y,z\in R$ one has $(xy,z)=(y,x^*z)$;
  \item For all $x\in R$ one has $\Vert x^*\Vert =\Vert x\Vert$;
  \item $xx^*\neq 0$ for $x\neq 0$.
\end{enumerate}
\end{defn}
An example of a Hilbert algebra is the set of Hilbert-Schmidt operators
$HS(I)$ that one can identify with the set of all matrices
$(a_{\alpha\beta})$, where $\alpha,\beta$ belong to a fixed set of
indices $ I$, satisfying the condition
$\sum_I|a_{\alpha\beta}|^2<\infty$.
\begin{thm} \cite{[N]}, p. 331. Every Hilbert algebra is the direct orthogonal
sum of its closed minimal two-sided ideals, which are simple Hilbert
algebras.

Every simple Hilbert algebra is isomorphic to some algebra $HS(I)$ of
Hilbert-Schmidt operators.
\end{thm}
\begin{defn}\cite{[L]}, p. 101
An idempotent $e\in R$ is said to be irreducible if it cannot be
expressed as a sum $e=e_1 + e_2$ with $e_1 ,\, e_2$ idempotents
which annihilate each other: $e_1e_2=e_2e_1=0$.
\end{defn}
\begin{lem}\cite{[L]}, p. 102.
A subset $I$ of a Hilbert algebra $R$ is a minimal left (right)
ideal if and only if it is of the form $I=R\cdot e$ ($I=e\cdot R$),
where $e$ is an irreducible self-adjoint idempotent. Moreover
$e\cdot R\cdot e$ is isomorphic to the set of complex numbers and
$R$ is spanned by its minimal left ideals.
\end{lem}
\noindent Observe that any minimal left ideal is closed, since it is
of the form $R\cdot e$.
\begin{cor}
If $R$ is a commutative Hilbert algebra, then any minimal left (or right) ideal is one-dimensional.
\end{cor}
\subsection{An algebra structure on $L^2(X)$}
Let $G=SL(n,\mathbb R)$, $H=GL(n-1,\mathbb R)$, $K=SO(n)$ and
$M=SO(n-1)$. We consider $H$ as a subgroup of $G$, consisting of the
matrices of the form $\left(
                         \begin{array}{cc}
                         (\det h)^{-1} & 0 \\
                         0 &   h   \\
                         \end{array}\right)$ with $ h\in GL(n-1,\mathbb R).$

Let $P^-$ be the parabolic subgroup of $G$ consisting of $1\times
(n-1)$ lower block matrices $P=\left(
                               \begin{array}{cc}
                                 a & 0 \\
                                 c & A \\
                               \end{array}
                             \right),\, a\in\mathbb R^*,\,c\in\mathbb R^{n-1}$ and
                             $A\in GL(n-1,\mathbb R)$ such that
                             $a\cdot\det A=1$.
 Similarly, let $P^+$ be the group of upper block matrices $P=\left(
                                                              \begin{array}{cc}
                                                                a & b \\
                                                                0 & A \\
                                                              \end{array}
                                                            \right)\, a\in\mathbb R^*,\,b\in\mathbb R^{n-1}$ and
                             $A\in GL(n-1,\mathbb R)$ such that
                             $a\cdot\det A=1$.

 The group $G$ acts on the sphere $S=\{s\in\mathbb R^n,\,\Vert
 s\Vert^2=1\}$ and acts transitively on the set $ \widetilde{S}=S/\sim$,
 where $s\sim s'$ if and only if $s=\pm s'$, by
 $\displaystyle g.s=\frac{g(s)}{\Vert g(s)\Vert}$, where $g(s)$ denotes the linear
 action of $G$ on $\mathbb R^n$.
 Clearly the stabilizer of the equivalence class of the first basis vector
 $\widetilde{e_1}$ is the group $P^+$, thus $\widetilde{S}\simeq
 G/P^+$. If $ds$ is the usual normalized surface measure on $S$,
 then $d(g.s)=\Vert g(s)\Vert^{-n}ds$.

  For $\mu\in\mathbb C$, define
 the character $\omega_\mu$ of $P^\pm$ by $\omega_\mu(P)=\vert
 a|^\mu$. Consider the induced representations
 $\pi_\mu^\pm={\mathrm{Ind}}_{P^\pm}^G\omega_{\mp\mu}$.

 Both $\pi_\mu^+$ and $\pi_\mu^-$ can be realized on
 $C^\infty(\widetilde{ S})$, the space of even smooth functions $\phi$
 on $S$. This action is given by
 $$
 \pi_\mu^+(g)\phi(s)=\phi(g^{-1}.s)\cdot\Vert g^{-1}(s)\Vert^\mu.
 $$
 Let $\theta$ be the Cartan involution of $G$ given by
  $\theta(g)={}^tg^{-1}$. Then
 $$
 \pi_\mu^-(g)\phi(s)=\phi(\theta(g^{-1}).s)\cdot\Vert\theta(g^{-1})(s)\Vert^\mu.
 $$

 Let $(\,,\,)$ denote the usual inner product on $L^2(S)$ :
 $(\phi,\psi)=\int_S\phi(s)\bar\psi(s)ds$. Then this sesqui-linear form
 is invariant with respect to the pairs of representations
 $(\pi_\mu^+,\pi_{-\mu-n}^+)$ and $(\pi_\mu^-,\pi_{-\mu-n}^-)$.
Therefore the representations $\pi_\mu^\pm$ are unitary for
${\mathrm{Re}}\,\mu=-\frac n2$.

The group $G$ acts also on $\widetilde{S}\times\widetilde{S}$ by
\begin{equation}\label{11}
g(u,v)=(g.u,\theta(g)v).
\end{equation}
This action is not transitive: the orbit
$(\widetilde{S}\times\widetilde{S})^{o}\,=\,G.(\widetilde{e_1},\widetilde{e_1})=\{(u,v)\,:\langle
u,v\rangle\neq 0\}/\sim$ is dense (here $\langle\,,\,\rangle$
denotes the canonical inner product on $\mathbb R^n$). Moreover
$(\widetilde{S}\times\widetilde{S})^o\simeq X$.\bigskip

The map $\displaystyle f\mapsto f(u,v)|\langle u,v\rangle|^{-\frac
n2+i\sigma}$, with $\sigma\in\mathbb R$ is a unitary $G$-isomorphism
between $L^2(X)$ and $\pi_{-\frac
n2+i\sigma}^+\,\widehat{\oplus}_2\,\pi_{-\frac n2+i\sigma}^-$ acting
on $L^2(\widetilde{S}\times\widetilde{S})$. The latter space is
provided with the usual inner product.

Define the operator $A_\mu$ on $C^{\infty}(\widetilde{S})$ by the
formula : $$ A_{\mu}\phi(s)=\int_S\vert\langle
s,t\rangle\vert^{-\mu-n}\phi(t)dt.$$
 This integral
converges absolutely for ${\mathrm{Re}}\,\mu<-1$ and can be
analytically extended to the whole complex plane as a meromorphic
function of $\mu$. It is easily checked that $A_{\mu}$ is an
intertwining operator, that is, $\displaystyle
A_{\mu}\pi_{\mu}^{\pm}(g)=\pi_{-\mu-n}^{\mp}(g)A_{\mu}$.

The operator $ A_{-\mu-n}\circ A_\mu$ intertwines the representation
$\pi_\mu^\pm$ with itself and is therefore a scalar
$c(\mu){\mathrm{Id}}$ depending only on $\mu$. It can be computed
using $K$-types.

Let $e(\mu)=\int_S|\langle s,t\rangle|^{-\mu-n}dt$, then
$c(\mu)=e(\mu)e(-\mu-n)$. But on the other hand side
$e(\mu)=\frac{\Gamma\left(\frac
n2\right)}{\sqrt\pi}\frac{\Gamma\left(\frac{-\mu-n+1}2\right)}{\Gamma\left(-\frac
{\mu}2\right)}$. One also shows that $A_\mu^*=A_{\bar\mu}$. So that, for
$\mu=-\frac n2+i\sigma$ we get (by abuse of notations):
$$
c(\sigma)=\frac{\Gamma\left(\frac
n2\right)^2}{\pi}\cdot\frac{\Gamma\left(\frac{-n/2-i\sigma
+1}2\right)\Gamma\left(\frac{-n/2+i\sigma +1}2\right)}
{\Gamma\left(\frac{n/2+i\sigma}2\right)\Gamma\left(\frac{-n/2-i\sigma}2\right)},
$$
and moreover $A_{-\frac n2+i\sigma}\circ A^*_{-\frac
n2+i\sigma}=c(\sigma ){\mathrm{Id}}$, so that the operator $d(\sigma )A_{-\frac
n2+i\sigma}$, where $d(\sigma )=\frac{\sqrt\pi}{\Gamma\left(\frac
n2\right)}\frac{\Gamma\left(\frac{n/2+i\sigma}2\right)}{\Gamma\left(\frac{-n/2+i\sigma +1}2\right)}$
is a unitary intertwiner between $\pi_{-\frac n2+i\sigma}^-$ and
$\pi_{-\frac n2-i\sigma}^+$.

We thus get a $\pi_{-\frac
n2+i\sigma}^+\,\widehat{\oplus}_2\,\bar\pi_{-\frac n2+i\sigma}^+$ invariant
map from $L^2(X)$ onto $L^2(\widetilde{S}\times\widetilde{S})$ given
by
$$
f\mapsto d(\sigma )\int_S f(u,w)|\langle u,w\rangle|^{-\frac
n2+i\sigma}|\langle v,w\rangle|^{-\frac n2-i\sigma}dw=:(T_\sigma f)(u,v),\,\forall
\sigma\neq0.
$$
This integral does not converge absolutely, it must be considered as
obtained by analytic continuation.

\begin{defn} A symbolic calculus on $X$ is a linear map $Op_\sigma : L^2(X)\to \mathcal{L}(L^2(\widetilde
S))$ such that for every $f\in L^2(X)$ the function $(T_\sigma f)(u,v)$ is
the kernel of the Hilbert-Schmidt operator $Op_\sigma (f)$ acting on
$L^2(\widetilde S)$.
\end{defn}
\begin{defn} The product $\#_\sigma$ on $L^2(X)$ is defined by
$$Op_\sigma (f\sharp_\sigma g)=Op_\sigma (f)\circ Op_\sigma (g),\,\forall f,g\in L^2(X).$$
\end{defn}
We thus have\begin{itemize}
              \item The product $\sharp_\sigma$ is associative.
              \item $\Vert f\sharp_\sigma g\Vert_2\leq\Vert f\Vert_2\cdot\Vert g\Vert_2$, for all $f,g\in L^2(X)$.
              \item $Op_\sigma (L_x f)=\pi_{-\frac n2+i\sigma}^+(x)\,
              Op_\sigma (f)\,\pi_{-\frac n2+i\sigma}^+(x^{-1}),$ so
              $L_x(f\sharp_\sigma\, g)=(L_xf)\sharp_\sigma\,(L_xg)$, for all
              $x\in G$, where $L_x$ denotes the left translation by
              $x\in G$ on $L^2(X)$.
            \end{itemize}

This non-commutative product can be described explicitly:
\begin{equation}\label{12}
    (f\sharp_\sigma\,g)(u,v)=d(\sigma )\int_S\int_S f(u,x)g(y,v)|[u,y,x,v]|^{-\frac
    n2+i\sigma}d\mu(x,y),
\end{equation}
where $d\mu(x,y)=|\langle x,y\rangle|^{-n}dxdy$ is a $G$-invariant
measure on $\widetilde S\times\widetilde S$ for the $G$-action
(\ref{11}), and $\displaystyle[u,y,x,v]=\frac{\langle
u,x\rangle\langle y,v\rangle}{\langle u,v\rangle\langle
x,y\rangle}.$

On the space $L^2(X)$ there exists an (family of) involution
 $f\to f^*$ given by : $Op_\sigma
(f^*)=:Op_\sigma (f)^*$. Notice that the correspondance $f\to
Op_\sigma (\bar f^*)$ is what one calls in pseudo-differential
analysis "anti-standard symbolic calculus". The link between symbols
of standard and anti-standard calculus in the setting of the
para-Hermitian symmetric space $X$ has been made explicit in
\cite{[PU]} Corollary 1.4, see also Section 3.

Obviously we have $(f\,\sharp_\sigma  g)^*=g^*\sharp_\sigma f^*$ and
with the above product and involution, the Hilbert space $L^2(X)$
becomes a Hilbert algebra.\bigskip

\section{The structure of the subalgebra of $K$-invariant functions in $L^2(X)$}
Let $\mathcal{A}$ be the subspace of all $K$-invariant functions in
$L^2(X)$.
\begin{thm}
The subset $\mathcal A$ is a closed subalgebra of $L^2(X)$ with respect to
the product $\sharp_\sigma$.
\end{thm}
This statement clearly follows from the covariance of the symbolic
calculus $Op_\sigma$, namely: $L_x(f\sharp_\sigma\,
g)=(L_xf)\sharp_\sigma\,(L_xg)$, for all
              $x\in G, f,g\in L^2(X)$.
\begin{thm}
Let $n>2$, then the subalgebra $\mathcal A$ is commutative.
\end{thm}
\emph{Proof.} For a function $f\in L^2(X)$ we set
$\check{f}(u,v)=f(v,u)$. The map $f\to\check{f}$ is a linear
involution. Indeed,
$$
(f\,\sharp_\sigma
g)(u,v)=d(\sigma )\int_{S}\int_S\check{f}(x,u)\check{g}(v,y)|[u,y,x,v]|^{-\frac
    n2+i\sigma}d\mu(x,y).$$
    Permuting $x$ and $y$ and $u$ and $v$ respectively, we get
$$
(f\,\sharp_\sigma g)(v,u)=d(\sigma )\int_{S
}\int_S\check{g}(u,x)\check{f}(y,v)|[v,x,y,u]|^{-\frac
    n2+i\sigma}d\mu(x,y).$$
    But $|[v,x,y,u]|=|[u,y,x,v]|$, therefore $(f\,\sharp_\sigma
    g)\check{}=\check g\,\sharp_\sigma\check f.$

    On the other hand, given a couple $(u,v)\in \widetilde{S}
    \times\widetilde{S}$ there exists an element $k\in K$ such that
    $k.(u,v)=(v,u)$. Geometrically $k$ can be seen as a rotation of
    angle $\pi[2\pi]$ around the axis defined by the bisectrix of
    vectors $u$ and $v$ in the plane they generate. Of course, such
    a $k$ exists for an arbitrary couple $(u,v)$ only if $n>2$.

    Hence for every $f\in\mathcal A$ we have $f=\check f$ and
    therefore $f\,\sharp_\sigma g=g\,\sharp_\sigma  f$, for $f,g\in\mathcal
    A.\quad\Box$
    \section{ Irreducible self-adjoint idempotents of $\mathcal A$}
We begin with a {\bf reduction theorem} for the multiplication and involution in $L^2(X)$.

As usual, we shall identify $L^2(X)$ with $L^2(\widetilde
S\times\widetilde S ;|\langle x,y\rangle|^{-n}dxdy )$. If $\phi\in
L^2(X)$ we shall write $\phi (u,v)=|\langle u,v\rangle |^{n/2 -i\sigma}\,
\phi_o(u,v)$. Then $\phi_o\in L^2(\widetilde S\times\widetilde S
;dsdt) =L^2(\widetilde S\times\widetilde S )$, and therefore the map
$\phi\to \phi_o$ is an isomorphism.
\begin{thm}
Under the isomorphism $\phi\to\phi_o$ the product $\#_\sigma$ translates
into
$$\phi_o\#_\sigma \psi_o\, (u,v)= d(\sigma )\int_S\int_S \phi_o(u,x)\, \psi_o (y,v) \, |\langle x,y\rangle |^{-n/2 -i\sigma}dxdy$$
and the involution becomes:
$$\phi_0^\ast (u,v)= \overline{d(\sigma )}^2\int_S\int_S \bar\phi_0 (x,y) \left ( |\langle x,v\rangle ||\langle u,y\rangle |\right )^{-n/2 +i\sigma }dxdy.$$
\end{thm}
The proof is straightforward. So we have translated the algebra
structure of $L^2(X)$ to $L^2(\widetilde S\times\widetilde S)$.

Let $\phi$ be an irreducible self-adjoint idempotent in $\mathcal
A$. We shall give an explicit formula for the $\phi_o$-component of
$\phi$.

Consider the decomposition of the space
$L^2(\widetilde{S})=\oplus_{\ell\in 2\mathbb N}V_\ell$, where
$V_\ell$ is the space of harmonic polynomials on $\mathbb R^n$,
homogeneous of even degree $\ell$.

Then the space $L^2(\widetilde{S}\times\widetilde{S})$ decomposes
into a direct sum of tensor products $\displaystyle\oplus_{\ell,m\in
2\mathbb N}V_\ell\otimes\bar V_m$ and consequently
$L^2_K(\widetilde{S}\times\widetilde{S})=\oplus_{\ell\in 2\mathbb
N}(V_\ell\otimes\bar V_\ell)^K$, where the sub(super-)script $K$ means: ``the $K$-invariants in''.

Let $\dim V_\ell=d$ and $f_1,\ldots,f_d$ be an orthonormal basis of
$V_\ell$. Then the function $\theta_\ell(u,v)=\sum_{i=1}^d
f_i(u)\bar f_i(v)$, that is the reproducing kernel of $V_\ell$, is,
up to a constant, the $K$-invariant element of
$V_\ell\otimes\bar V_\ell$.

\begin{thm}
Let $\phi (u,\, v)=|\langle u,\, v\rangle |^{n/2 -i\sigma}\phi_o (u,\, v)$ be an irreducible self-adjoint idempotent in $\mathcal A$. Then there exist complex numbers $c(\sigma ,\, \ell)$ such that for any $\ell\in 2\mathbb N$ one has
$$\phi_o (u,\, v)=c(\sigma ,\, \ell )\, \theta_\ell (u,\, v).$$
For different $\ell$ and $\ell'$ the idempotents annihilate each other. Moreover they span $\mathcal A$.
\end{thm}
\emph{Proof.} Firstly we shall show that $\theta_\ell$ satisfies the condition
$$\theta_\ell \#_\sigma \theta_\ell =a(\sigma,\,\ell )\,\theta_\ell$$
for some constant $a(\sigma ,\, \ell )$. Indeed,
\begin{eqnarray*}
&\ &d(\sigma )\,\int_S\int_S \theta_\ell (u,x)\,\theta_\ell (y,v)\, |\langle x,\, y\rangle|^{-\frac n2 -i\sigma} dxdy\\
&=& d(\sigma )\, e_\ell (\sigma )\int_S \theta_\ell (u,y)\,\theta_\ell (y,v)\, dy = d(\sigma )\, e_\ell (\sigma )\,\theta_\ell (u,v)
\end{eqnarray*}
by the intertwining relation (apply $A_{-\frac n2 +i\sigma}$ to $\theta_\ell (.,x)$):
$$\int_S \theta_\ell (u,x)\, |\langle x,\, y\rangle |^{-\frac n2 -i\sigma} dx=e_\ell (\sigma )\,\theta_\ell (u,y)$$
where $e_\ell (\sigma )=\int_S \displaystyle\frac{\theta_\ell (e_1,\, x)}{\theta_\ell (e_1,\, e_1)}\, |x_1|^{-\frac n2 -i\sigma}\, dx.$

Observe that $\displaystyle\frac{\theta_\ell (e_1,\, x)}{\theta_\ell (e_1,\, e_1)}$ is a spherical function on $\widetilde S$ with respect to $M$ of the form $a_\ell\, C_\ell^\frac{n-2}2 (|x_1|)$ where $C_\ell^\frac{n-2}2(u)$ is a Gegenbauer polynomial and
$$a_\ell^{-1} =C_\ell^\frac{n-2}2 (1)=2^\ell \displaystyle \frac{\Gamma (\frac{n-2}2 +\ell )}{\Gamma (\frac{n-2}2 )\ell!}.$$
See for instance \cite{[V]}, Chapter IX, \S 3. Notice that $\theta_\ell (e_1,\, e_1)={\rm dim}\, V_\ell =\displaystyle\frac {(n+\ell -1)!}{(n-1)!\ell!}\neq 0$. The integral defining $e_\ell (\sigma )$ does not converge absolutely, but has to be considered as the meromorphic extension of an analytic function. Poles only occur in half-integer points on the real axis. So we have to restrict (and we do) to $\sigma\neq 0$.

So we have $\theta_\ell \#_\sigma \theta_\ell = d(\sigma )\, e_\ell (\sigma )\,\theta_\ell$ and hence $\varphi_\ell = [d(\sigma )\, e_\ell (\sigma )]^{-1}\,\theta_\ell$ is the $\phi_o$-component of an idempotent in $\mathcal A$. Furthermore $\theta_\ell \#_\sigma \theta_{\ell'}=0$ if $\ell\neq\ell'$. Clearly $\varphi_\ell$ is self-adjoint, since $|d(\sigma )|^{-2}=|e_\ell (\sigma )|^2$, being equal to the constant $c(\sigma )$ from Section 1.

So the $\varphi_\ell$ are mutually orthogonal idempotents in the
algebra $L_K^2((\widetilde S\times\widetilde S ); dsdt)$ and span
this space. The theorem now follows easily. $\quad\Box$
\vskip.2cm\noindent {\bf Remark} The constant $e_\ell (\sigma )$ can
of course be computed. Apply e.g. \cite{[G]}, Section 7.31, we get,
by meromorphic continuation:
\begin{eqnarray*}
e_\ell (\sigma )&=&a_\ell \, \int_S C_\ell^{\frac {n-2}2}(|x_1|)\, |x_1|^{-\frac n2 -i\sigma}\, dx\\
&=& 2\, a_\ell\, \frac{\Gamma ({\frac n2})}{\Gamma ({\frac {n-1}2})\, \sqrt\pi}\, \int_0^1 u^{-\frac n2 -i\sigma}\, (1-u^2)^{\frac {n-2}2}\, C_\ell^{\frac {n-2}2}(u)\, du\\
&=& 2^{-2\ell}\,\frac{\Gamma (\frac n2 )}{\sqrt\pi}\cdot \frac{\Gamma (n-2+\ell )}{\Gamma (n-2)}\cdot \frac{\Gamma (\frac{n-2}2)}{\Gamma (\frac{n-2}2 +\ell)}\cdot \frac{\Gamma (-\frac n2 -i\sigma +1)\Gamma (\frac{-\frac n2 -i\sigma -\ell +1}2)}{ \Gamma (-\frac n2 - i\sigma -\ell +1)\Gamma (\frac{\frac n2 -i\sigma +\ell }2)}.
\end{eqnarray*}

\section{The strucure of the Hilbert algebra $L^2(X)$}
We now turn to the full algebra $L^2(X)$. We again reduce the computations to $L^2(\widetilde S\times \widetilde S )$. In a similar way as for $\mathcal A$ we get:
\begin{lem}\label{lemma} If $\phi_o\in V_\ell\otimes \overline V_m,\, \psi_o\in V_{\ell'}\otimes\overline V_{m'}$ then
$$\phi_o \#_\sigma\, \psi_o = \left\{\begin{array}{ll} 0 &\mbox{{\rm if}\, $m\not= \ell'$}\\ {\rm in}& \mbox{$V_\ell\otimes \overline V_{m'}$ if $m=\ell'$}.
\end{array}\right.$$
More precisely we have the following result. Let $(f_i),\, (g_j),\,
(k_l)$ be orthonormal bases of $V_\ell,\, V_m$ and $V_{m'}$
respectively, and $\phi_o (u,v)=f_i(u)\overline g_j(v),\ \psi_o
(u,v)=g_{j'}(u)\overline k_l(v)$, then
$$\phi_o \#_\sigma\, \psi_o = \left\{\begin{array}{ll} 0 &\mbox{if $j\not= j'$}\\  d(\sigma )\, e_m(\sigma )\, f_i(u)\overline k_l(v)
& \mbox{if $j=j'$}.\end{array}\right.$$
\end{lem}
The proof is again straightforward and uses the intertwining
relation:
$$\int_S |\langle x,y\rangle |^{-n/2-i\sigma}\, g_{j'}(y) dy = e_m(\sigma )\, g_{j'}(x).$$
\begin{thm}
The irreducible self-adjoint idempotents of $L^2(\widetilde
S\times\widetilde S)$ are given by
$$e^\ell_f(u,v)= \{d(\sigma )\, e_\ell(\sigma )\}^{-1}\cdot f(u)\, \overline f(v)$$
with $f\in V_\ell,\, \Vert f\Vert_{L^2(\widetilde S )}=1$ and $\ell$
even. The left ideal generated by $e^\ell_f$ is equal to
$L^2(\widetilde S)\otimes \overline f$.
\end{thm}
The proof is by application of Lemma (\ref{lemma})
\vskip.2cm\noindent {\bf Remarks} \begin{enumerate}
 \item The minimal
right ideals are obtained in a similar way. \item The minimal
two-sided ideal generated by $L^2(\widetilde S\times\widetilde
S)\cdot e^\ell_f$ is the full algebra $L^2(\widetilde
S\times\widetilde S )$. \item The closure of
$\bigoplus_{\ell\in2\mathbb N}V_\ell \otimes \overline V_\ell$ is a
$H^*$-subalgebra of $L^2(\widetilde S \times\widetilde S )$. The
minimal left ideals are here $V_\ell\otimes \overline f\ (f\in
V_\ell, \, \Vert f\Vert _{L^2(\widetilde S )}=1)$; they are
generated by the $e^\ell_f$ as above. The minimal two-sided ideal
generated by $V_\ell\otimes \overline f$ is equal to $V_\ell\otimes
\overline V_\ell$.
\end{enumerate}
\section{The case of a general para-hermitian space}
It is not necessary to assume ${\rm rank}\, X=1$ in order to show that $\mathcal A$ is commutative. Theorem 3.2
 is also valid mutatis mutandis in the general case since $(K,\, K\cap H)$ is a Gelfand pair, and it clearly implies the commutativity of $\mathcal A$. To the general construction of the product and the involution we shall return in another paper.
\bibliographystyle{amsplain}

\end{document}